\input AHTOH-E.STY
\hfuzz6.5pt

\def\frac#1#2{{{#1}\over{#2}}}

\UDC{512.543.7+512.542}
\title{Economical adjunction of square roots to groups}
\author{Dmitrii V. Baranov \quad Anton A. Klyachko}
\address{Faculty of Mechanics and Mathematics\\
       Moscow State University\\
       Moscow 119991, Leninskie gory, MSU\\
       dimbaranov@mail.ru\quad klyachko@mech.math.msu.su
}
\grants{\RFBR08-01-00573}

\abstract{%
\narrower
\narrower
\narrower
How large must an overgroup of a given group be in order to contain 
a square root of any element of the initial group? We give an 
almost exact answer to this question 
(the obtained estimate 
is at most twice worse than the best possible) 
and state several related open questions.  
}

\s 0.
Introduction

The solvability of equations over groups has been extensively
studied 
(see, e.g., [GeRo62], [Levin62], [Lynd80], [Brod84], 
[EdHo91], [Howie91], [Klya93], [KlPr95], [FeRo96], [Klya97], [ClGo00], 
[EdJu00], [Juh\'a03], [Klya06] and the references therein). In these 
papers, it was proven that, under some conditions, an equation $w(x)=1$ 
with coefficients from a group $G$ is \emph{solvable over $G$}, i.e., 
there exists a group $H$ containing $G$ as a subgroup and an 
element~${h\in H}$ such that $w(h)=1$. In this paper, 
we study a quantitative 
question: {\sl how large must such a group $H$ be?} Even for 
simple equations whose solvability is well known, this question 
turns out to be very difficult and we restrict ourselves to  
the simplest nontrivial equation~${x^2=g}$.

Certainly, the answer strongly depends on the initial group $G$. For 
example, if the order of $G$ is odd, then the role of $H$ can be played by 
the group $G$ itself; if $G$ is cyclic, then $H$ may be taken  
twice larger than $G$, etc. The most interesting is  
to estimate the order of $H$ ``in the worst case". 
We obtain an estimate at most twice worse than the best 
possible. Namely, we prove the following theorem.

\proclaim {Main theorem}.
Each finite group $G$ embeds into a group of order $2|G|^2$ in 
which all elements of $G$ are squares. There exist infinitely  
many pairwise nonisomorphic finite groups~$G_i$ and elements~$g_i \in G_i$ 
such that no group of order smaller than $|G_i|^2$ can 
contain $G_i$ together with a quadratic root of $g_i$.

Actually, we consider two problems: economical adjunction of a solution to 
an equation $x^2=g$ and economical adjunction of solutions to 
all equations of such form. The main theorem shows that, in both cases, 
a group of order $2|G|^2$ is enough, but a group of order less 
than $|G|^2$ is not enough sometimes.

The behaviour of solution sets of equations in finite groups 
has been throughly studied
(see, e.g., [Frob03], [Hall36], [Solo69], [Stru95], and 
the references therein). Unfortunately, we are unable to use 
these nontrivial results. 

The first assertion of the theorem is not 
new and can be proven easily (see~Section 1). In Section~2, we 
prove the second assertion.
In the last section, we state several open questions about 
economical adjunction of solutions of equations to groups.

\smallskip
\noindent
{\bf Notation} 
which we use is mainly standard. Note only 
that if $k\in \Z$ and $x$ and $y$ are elements of a group, then $x^y$, 
$x^{ky}$, and $x^{-y}$ denote $y^{-1}xy$, $y^{-1}x^ky$ and 
$y^{-1}x^{-1}y$, respectively.  If $X$ is a subset of a group, then 
$|X|$, $\gp X$, and $\nc X$ denote the cardinality of $X$, 
the subgroup generated by $X$, and the normal subgroup generated by 
$X$, respectively. The letter~$\Z$ denotes the set of integers. 
The symbol $\Z_n$ denotes the group or 
the ring $\Z/n\Z$. The multiplicative group of the ring $\Z_n$ 
is denoted by $\Z_n^*$.
The automorphism group of $G$ 
is denoted by $\Aut G$. The symbol~$D_p$ denotes the dihedral group of 
order $2p$. The stabiliser of a point $a$ under an action of a group $G$ 
is denoted by $\St_G(a)$.  \emph{A reflection} is an element of $D_p$ not 
belonging to its subgroup $\Z_p$.

The authors thank an anonymous referee for useful remarks.

\s 1.
Wreath products and the proof of the first assertion of the theorem

The first assertion of the theorem is well known [Levin62]:
the wreath product 
$$
G\wr\Z_2= 
\left\{
\pmatrix{
g_1&0\cr
0&g_2
}
\;\Big|\; 
g_1,g_2\in G 
\right\} 
\cup 
\left\{
\pmatrix{
0&g_1\cr
 g_2&0
}
\;\Big|\; 
g_1,g_2\in G 
\right\}
$$
of $G$ and the cyclic group of order~2
is a group of order $2|G|^2$ containing square roots of
all elements of $G$, assuming that group $G$ embeds into
the wreath product $G\wr\Z_2$ as the diagonal: 
$$
g\mapsto 
\pmatrix{
g&0\cr
0 &g
}.
$$ 
Indeed,
$
\left({{\scriptscriptstyle0}\;g \atop
{\scriptscriptstyle1}_{\;}{\scriptscriptstyle0}}\right)^2=
\left({g\;{\scriptscriptstyle0} \atop
{\scriptscriptstyle0}_{\;}g}\right). 
$
This is the simplest special case of the Levin theorem. The full  
statement can be found in the last section of this paper. 

\s 2. Dihedral groups and the proof of the second assertion of the theorem

The second assertion of the main theorem is implied by the following fact.

\Th 1.
If $p\in4\Z+3$ is prime, $\~G$ is a group containing the dihedral subgroup 
$G=D_p$, and the reflection $g\in G$ is the square of some element 
$x\in\~G$, then $|\~G|\ge |G|^2$.

\smallskip

To prove this theorem, we need several simple lemmata.

\Lemma 1. 
If $H_1$ and $H_2$ are subgroup of a group $H$, then
$|H| \ge {|H_1||H_2|\over|H_1\cap H_2|}=|H_1H_2|$.

We leave the proof of this lemma to the readers.

\Lemma 2.
If $D_p=G \subseteq \gp{G, x}=\~G$ and
$x^2=g$, where $g \in G$ is a reflection, then either 
$G \nin \~G$ or $G \cap G^x=\gp g$.

\Proof
Clearly, $g\in G \cap G^x$. The group $D_p$ has exactly 
two subgroups containing $g$. If $G \cap G^x=\gp g$, then
we have nothing to prove. If $G \cap G^x=G$, then $G=G^x$ and hence 
$G \nin \~G$, because $\gp{G, x}=\~G$.

\Lemma 3.
If $D_p=G \nin \~G$, where $p\in 3+4\Z$ is
prime, then no reflection $g \in G$ 
is a square in $\~G$.

\Proof
The subgroup $\Z_p\subset D_p=G \nin \~G$
is the commutator subgroup of $G$ and, therefore, it is characteristic in 
$G$ and normal in $\~G$.
The group $\~G$ acts on $\Z_p$ by conjugations. 
The reflection $g$ acts as 
$-1\in\Z_p^*=\Aut\Z_p$, and it is well known that $-1$ is not a  
square in $\Z_p^*$ if $p\in 3+4\Z$. This completes the proof.

\bigskip

Now, we proceed to prove the main theorem. We may 
assume that $\~G=\gp{G,x}$. Let
$K$ be the set of all subgroups of $\~G$ conjugate with $G$. 
The group $\~G$ acts transitively on $K$ by conjugations.

\Lemma 4.
$|\~G|\ge |K|\cdot|G|$.

\Proof
$|\~G|=|K|\cdot|\St_{\~G}(G)|\ge |K|\cdot|G|$, because 
$G \subseteq\St_{\~G}(G)$.
\bigskip

Consider the complete 
undirected graph $\Gamma$ with vertex set $K$. 
We call an edge $(G^{h_1}, G^{h_2})$
$$ 
\hbox{
green if $|G^{h_1}\cap G^{h_2}|=2$;
\quad 
yellow if $|G^{h_1}\cap G^{h_2}|=p$;
\quad 
red if $|G^{h_1}\cap G^{h_2}|=1$.
}
$$
Clearly, all edges are coloured.
If there is at least one red edge, then the assertion 
follows readily from Lemma 1. So, we assume that there are no
red edges. 

\Lemma 5.
All vertices $K$ and yellow edges form a graph $Y$ 
whose  
connected components are complete graphs. All 
these components 
contain the same number of vertices.

\Proof
The first assertion follows immediately from that  
$D_p$ has a unique subgroup of order $p$.
The second assertion 
follows from that the action of $\~G$ on $K$ is transitive and preserves 
colours of edges.

\Lemma 6.
The number of green edges incident to the vertex $G$ is a positive 
multiple of $p$.

\Proof
Each green edge incident to $G$ corresponds to one of the $p$ 
reflections $g\in G$. Thus, the edges incident to $G$ are divided into $p$ 
classes. These classes have the same number of edges, because 
all reflections in $G$ are conjugate and, therefore, 
any class is mapped onto any other class by
an automorphism of 
the graph. 
This means that the number of green edges incident to $G$ is a multiple 
of $p$.

The graph has a green edge, namely, the edge $(G,G^x)$, 
where $x^2=g\in G$ is a reflection. Indeed, $G^x\cap G\ni g$, 
but $G^x\ne G$ (because otherwise $G$ would be a normal subgroup in 
$\~G=\gp{G,x}$, which contradicts Lemma 3). Therefore, 
$G^x\cap G= \gp{g}_2$ and the lemma is proven.    

\bigskip
We continue the proof of the theorem. 
Suppose that there are at least $2p$ green edges incident to $G$.
Then the graph has at least $2p+1$ vertices, i.e., $|K|>2p$, and, 
by Lemma 4, $|\~G|\ge |K|\cdot|G|>2p|G|=|G|^2$, as required.
 
According to Lemma 6, it remains to consider the case, where 
there are exactly $p$ green edges incident to each vertex.

Suppose that $u$ is the number of vertices in 
each connected component of $Y$ (Lemma 5), and
$v$ is the number of these components. 
Then 
$$
p=\hbox{(the number of green edges incident to $G$)}=(v-1)u 
$$
(because each vertex not joined with 
$G$ by a yellow edge is joined with $G$ by a green edge). 

The equality $p=(v-1)u$ means (by virtue of the primeness of $p$) that 
either $v=2$ and $u=p$ or $v=p+1$ and $u=1$.

In the first case, $|K|=2p$, and the assertion follows from Lemma 4:  
$|\~G|\ge 2p|G|=|G|^2$.

In the second case, $|K|=p+1$, and the graph $\Gamma$ is a complete graph 
all of whose edges are green. 
The group~$\~G$ acts on this graph, and the action of $G$ on 
set vertices different from $G$ is isomorphic to the action of $G$ 
by conjugations on the set of its subgroups of order two
(this isomorphism maps a group~$G^h$ to the subgroup $G^h\cap G$). 
In particular, the conjugation by a reflection $g$ 
is a permutation of vertices of the graph which fixes exactly two points
($G$ and $G^h$, where $G\cap G^h=\gp g$) 
and, therefore, this permutation decomposes into a product 
of $p-1\over2$ independent transpositions. This permutation is 
odd, because $p\in 3+4\Z$, which contradicts the assumption that $g$ is 
a square in $\~G$. The main theorem is proven.

\s 3.
Higher degree roots and other open questions

The question arises: what is the best possible estimate?

\Question 1. 
Do there exist infinitely many finite groups $G$ such that, 
for some $g\in G$, each overgroup of $G$ in which~$g$ is a square 
has order at least $2|G|^2$?

The following proposition shows that, for dihedral groups, our
theorem can not be strengthened, and
to answer Question~1 one must study groups close to simple.

\Proposition 1. If a finite group $G$ and its element $g$
satisfy at least one of the conditions 
\item{\rm a)}
$G$ does not coincide with its commutator subgroup; 
\item{\rm b)} 
$G$ does not coincide with the normal closure of $g$; 
\item{\rm c)}
$G$ has a nontrivial normal subgroup of odd order,%
\fn{%
\rm The Feit--Thompson theorem about the solvability of odd-order groups 
[FeTh63] implies that property c) is equivalent to the existence 
of a nontrivial abelian normal odd order subgroup of $G$.
}
\enditem 
then $G$
embeds into a group $H$ of order at most $|G|^2$ in
which the element $g$ is a square.

\Proof
The following lemma shows that under condition~a) or~b) a proper subgroup 
of the wreath product $G\wr\Z_2$ (see Section~1) can be taken for $H$. If 
condition~c) holds, then we can take a proper quotient of this wreath 
product for~$H$; this follows from Lemma 8 (see below).

\Lemma 7.
In the wreath product $G\wr\Z_2$, the subgroup $H$ generated by
$G$ embedded diagonally and the quadratic root
$\left({{\scriptscriptstyle0}\;
g\atop{\scriptscriptstyle1}_{\;}{\scriptscriptstyle0}}\right)$
of an element $g\in G$ has the form
$$
H= \left\{
\pmatrix{
g_1&0\cr
0&g_2}
\;;\;
g_1g_2^{-1}\in \[\nc g,G\] \right\} \cup \left\{
\pmatrix{
0&g_1\cr
g_2&0}
 \;;\; g_1g_2^{-1}\in g\[\nc g,G\] \right\},
$$
where $\[\nc g,G\]$ is the mutual commutator subgroup of 
$G$ and 
the normal closure
of $g$ in $G$.

\Proof
The natural epimorphism $\phi\:G\to G/\[\nc g,G\]$ induces a homomorphism 
$\Phi\:G\wr\Z_2\to \(G/\[\nc g,G\]\)\wr\Z_2$. The  
right-hand side of the required equality is $\Phi^{-1}(\Phi(H))$. 
Therefore, it is sufficient to prove that $H$ contains the kernel of
$\Phi$. But $\ker\Phi$ is generated (as a subgroup) by the elements 
of the form 
$$ 
\pmatrix{ \[g^x,y\]&0\cr 0&1 } \qqbox{and} \pmatrix{ 1&0\cr 
0&\[g^x,y\]
},
\qbox{where $x,y\in G$},
$$
which lie in $H$, as the following equalities shows:
$$
\eqalign{
\pmatrix{
x^{-1}&0\cr
0&x^{-1}
}
\pmatrix{
0&g\cr
1&0
}
\pmatrix{
x&0\cr
0&x
}
=
\pmatrix{
0&g^x\cr
1&0
},
\quad
\[
\pmatrix{
0&g^x\cr
1&0
}
,
\pmatrix{
y&0\cr
0&y
}
\]
=
\pmatrix{
0&1\cr
g^{-x}&0
}
\pmatrix{
0&g^{xy}\cr
1&0
}
&=
\pmatrix{
1&0\cr
0&\[g^x,y\]
},
\cr
\pmatrix{
1&0\cr
0&\[g^x,y\]
}
\pmatrix{
\[g^x,y\]&0\cr
0&\[g^x,y\]
}^{-1}
&=
\pmatrix{
\[g^x,y\]&0\cr
0&1
}^{-1}.
}
$$

\Lemma 8.
If $N$ is a normal abelian subgroup of a group $G$,
then the set
$$
K= \left\{
\pmatrix{
x&0\cr
0&x^{-1}
}
 \;;\; x \in N \right\}
$$
is a normal subgroup of the wreath product $G\wr\Z_2$.
If the order of $N$ is odd, then
the intersection of $K$
and $G$ (embedded into the wreath product diagonally) is trivial.
\newline
Conversely: each nontrivial normal subgroup of this
wreath product trivially intersecting $G$ contains a nontrivial
abelian normal subgroup of the specified form.

\Proof
The set $K$ is a normal subgroup, obviously. 
Clearly, $K\cap G=\{x\in N\;;\;x^2=1\}$; therefore, $K$ intersects $G$
trivially, if the order of $N$ is odd.

It is well known that 
an arbitrary nontrivial normal subgroup $X$ of a wreath product
nontrivially intersects the base (see, e.g., [KaMe82]). If 
$$
1\ne u=
\pmatrix{
x&0\cr
0&y
}
\in X,
$$
then
$$
\[u,
\pmatrix{
y&0\cr
0&y
}
\]=
\pmatrix{
[x,y]&0\cr
0&1
}
=v
\in X\ni 
\pmatrix{
0&1\cr
1&0
}^{-1}
v
\pmatrix{
0&1\cr
1&0
}
=
\pmatrix{
1&0\cr
0&[x,y]
}
=w
$$
and, therefore,
$$
vw=
\pmatrix{
[x,y]&0\cr
0&[x,y]
}
\in X\cap G=\1,
\qbox{i.e., $[x,y]=1$.}
$$
But then
$$
X\ni 
\pmatrix{
0&1\cr
1&0
}^{-1}
u
\pmatrix{
0&1\cr
1&0
}
=
\pmatrix{
y&0\cr
0&x
}
=t
$$
and, therefore,
$$
ut=
\pmatrix{
xy&0\cr
0&xy
}
\in X\cap G=\1,
\qbox{i.e., $xy=1$.}
$$
Thus, the intersection of $X$ and the base of the wreath product has  
the form 
$$
K= \left\{
\pmatrix{
x&0\cr
0&x^{-1}
}
 \;;\; x \in N \right\},
$$
where $N$ is a subset of $G$.
This implies that the set $N$ must be an abelian 
normal subgroup. Lemma 8 is proven.

\medskip
 
These lemmata prove Proposition 1 and, in addition, show that, if 
$G$ satisfies neither condition a), b), nor c) (e.g., 
if the group $G$ is nonabelian simple), then the wreath product $G\wr\Z_2$ 
has neither proper subgroups nor proper quotient 
containing $G$ and the square 
root~$\left({{\scriptscriptstyle0}\;g\atop{\scriptscriptstyle1}_{\;}{\scriptscriptstyle0}}\right)$
of~$g$. 

\medskip

Now, consider higher degree roots and solutions to other equations.
The point of departure in our investigation is the Levin theorem, whose 
complete statement looks as follows.

\proclaim{Levin theorem {\rm ([Levin62])}}. 
The wreath product $G\wr\Z_n$
\(whose order 
is~$n|G|^n$\) 
of a group $G$ and the cyclic group of order $n$ 
contains solutions to all positive equations of degree $n$ 
over the group $G$.

A \emph{positive equation of degree $n$ over a group $G$}
is an equation of the form
$$
g_1xg_2x\dots g_nx=1
\qbox{where $g_1,\dots,g_n\in G$}.
$$
So, the question arises: does the Levin theorem give the
best possible estimate?

\Question 2.
Do there exist
infinitely many finite groups $G$ such
that each overgroup of $G$ 
containing solutions to all positive equations of 
degree $n$ over $G$ has order at least $n|G|^n$?

One can formulate a more daring conjecture.

\Question 3.
Do there exist 
infinitely many finite groups $G$ such
that each overgroup $H$ of $G$ in which all elements of 
$G$ are $n$th powers \(of elements of $H$\) has order 
at least $n|G|^n$?

What can be said about economical adjunction of solutions to other 
(i.e., nonpositive) equations? For example, 
the Gerstenhaber--Rothaus theorem [GeRo62] in combination with 
Mal'tsev's theorem about 
residual finiteness of finitely 
generated linear groups [Mal40] gives the following fact.

\Proposition 2.
{\rm([GeRo62]+[Mal40])}.
Each finite group $G$
embeds into a finite group~$H$ containing solutions to all
nonsingular equations of length $n$ over $G$.

A \emph{nonsingular equation of length $n$ over a group $G$}
is an equation of the form
$$
g_1x^{\epsilon_1}g_2x^{\epsilon_2}\dots g_nx^{\epsilon_n}=1,
\qbox{where $g_i\in G$,\quad $\epsilon_i\in\{\pm1\}$,
\qbox{and }
$\sum\epsilon_i\ne0$}.
$$

The proof of the Gerstenhaber--Rothaus theorem is nice but 
nonconstructive. So, it is difficult to write out  
any estimate of the order of $H$.

\Question 4.
How to estimate $|H|$ via $|G|$ and $n$ in
Proposition 2?

For $n=1$, the answers to Questions 2, 3, and 4 is, obviously, positive.
A nonsingular equation of length two has the form 
$g_1x^\epsilon g_2x^\epsilon=1$, where $\epsilon\in\{\pm1\}$, and a linear 
change of variables reduces it to the form $x^2=g$; therefore, the main 
result of this paper gives an answer to the ``twice weakened" versions of 
these questions for $n=2$. What occurs for other $n$ we do not know.

\REFERENCES

\[Brod84]
Brodskii S.D.
Equations over groups and one-relator groups
{// Sib. Mat. Zh.} 1984. {T.25}. no.2. P.84--103.

\[KaMe82]
Kargapolov M.I., Merzlyakov Yu.I.
Fundamentals of group theory. 
Moscow: ``Nauka", 1982.

\[Klya06a]
Klyachko Ant. A.
How to generalize known results on equations over groups
{// Mat. Zametki}. 2006. {T.79}. no.3. P.409--419.
See also
arXiv:math.GR/0406382.

\[KP95]
Klyachko Ant. A., Prishchepov M. I.
The descent method for equations over  groups
{// Moscow Univ. Math. Bull.} 1995, {V.50}  P. 56--58.

\[Mal40]
Mal'tsev A.I.
On isomorphic representation of infinite groups by
matrices {// Mat. Sb.} 1940. 8(50):3  P.405--422.

\[Stru95]
Strunkov S.P.
To the theory of equations on finite groups
// Izv. Ross. Akad. Nauk. Ser. Mat. 1995.
V.59:6. P.171--180

\[ClGo00]
Clifford A., Goldstein R.Z.
Equations with
torsion-free coefficients {// Proc. Edinburgh Math. Soc.} 2000.
{V.43}. P.295--307.

\[EdHo91]
Edjvet M., Howie J.
The solution of length four equations over
groups {// Trans. Amer. Math. Soc.} 1991. {V.326}. P.345--369.

\[EdJu00]
Edjvet M., Juh\'asz A.
Equations of length 4 and one-relator
products {// Math. Proc. Cambridge Phil. Soc.} 2000 V.129.
P.217--230

\[FeTh63]
Feit W., Thompson J.G.
Solvability of groups of odd order {//
Pacific J. Math.} 1963. V.13 P.755--1029.

\[FeRo96]
Fenn R., Rourke C.
Klyachko's methods and the solution of
equations over torsion-free groups {//L'Enseignment
Math\'ematique.} 1996. {T.42}. P.49--74.

\[Frob03]
Frobenius G.
\"Uber einen Fundamentalsatz der Gruppentheorie
// Berl. Sitz. 1903. S.987--991.


\[GeRo62]
Gerstenhaber M., Rothaus O.S.
The solution of sets of equations in
groups {// Proc. Nat. Acad. Sci. USA}. 1962. {V.48} P.1531--1533.

\[Hall36]
Hall Ph.
On a Theorem of Frobenius
// Proc. London Math. Soc. 1936. V.40. P.468--531.

\[Howie91]
Howie J.
The quotient of a free product of groups by a single
high-powered relator. III: The word problem {// Proc. Lond. Math.
Soc.} 1991. {V.62}. No.3 P.590-606.

\[Juh\'a03]
Juh\'asz A.
On the solvability of a class of equations over
groups. {// Math. Proc. Cambridge Phil. Soc.} 2003. V.135.
P.211--217

\[Klya93]
Klyachko Ant.A.
A funny property of a sphere and equations over
groups {// Comm. Algebra}. 1993. {V.21}. P.2555--2575.

\[Klya97]
Klyachko Ant.A.
Asphericity tests {// J. Algebra}. 1997. {V.7}.
P.415--431.

\[Levin62]
Levin F.
Solutions of equations over groups {// Bull. Amer. Math.
Soc.} 1962. {V.68}. P.603--604.

\[Lynd80]
Lyndon R.C.
Equations in groups {// Bol. Soc. Bras. Math.} 1980.
{V.11}. no.1. P.79--102.

\[Solo69]
Solomon L. 
The solutions of equations in groups
// Arch. Math. 1969. V.20. no.3. P. 241--247.

\end